\documentclass[12pt,leqno]{amsart}

\usepackage{color}

\newtheorem{theorem}{Theorem}
\newtheorem{proposition}[theorem]{Proposition}

\newtheorem{remark}{Remark}
\newtheorem{lemma}[theorem]{Lemma}
\newfont{\bb}{msbm10 at 12pt}

\def\pf{\noindent{\textit {Proof.} }}

\def\R{\hbox{\bb R}}

\def\S{\hbox{\bb S}}

\newcommand{\bal}{\begin{align}}      \newcommand{\eal}{\end{align}}
\newcommand{\ba}{\begin{array}}      \newcommand{\ea}{\end{array}}
\newcommand{\bc}{\begin{center}}     \newcommand{\ec}{\end{center}}
\newcommand{\be}{\begin{enumerate}}  \newcommand{\ee}{\end{enumerate}}
\newcommand{\beq}{\begin{eqnarray}}  \newcommand{\eeq}{\end{eqnarray}}
\newcommand{\beQ}{\begin{eqnarray*}} \newcommand{\eeQ}{\end{eqnarray*}}
\newcommand{\bi}{\begin{itemize}}    \newcommand{\ei}{\end{itemize}}
\newcommand{\bt}{\begin{tabular}}    \newcommand{\et}{\end{tabular}}
\newcommand{\bdm}{\begin{displaymath}} \newcommand{\edm}{\end{displaymath}}

\def\qed{\hfill{q.e.d.}\smallskip\smallskip}

\begin{document}

\title[Uniqueness of the de Sitter spacetime]{Uniqueness of the de Sitter spacetime among static vacua with positive cosmological constant}

\author{Oussama Hijazi}
\address[Hijazi]{Institut {\'E}lie Cartan de Lorraine,\\
Universit{\'e} de Lorraine, Nancy \\
B.P. 239\\
54506 Vand\oe uvre-L{\`e}s-Nancy Cedex, France}
\email{Oussama.Hijazi@univ-lorraine.fr}

\author{Sebasti{\'a}n Montiel}
\address[Montiel]{Departamento de Geometr{\'\i}a y Topolog{\'\i}a\\
Universidad de Granada\\
18071 Granada \\
Spain}
\email{smontiel@ugr.es}

\author{Simon Raulot}
\address[Raulot]{Laboratoire de Math\'ematiques R. Salem
UMR $6085$ CNRS-Universit\'e de Rouen
Avenue de l'Universit\'e, BP.$12$
Technop\^ole du Madrillet
$76801$ Saint-\'Etienne-du-Rouvray, France.}
\email{simon.raulot@univ-rouen.fr}

\begin{abstract}
We prove that,  among all $(n+1)$-dimensional spin static vacua with positive cosmological constant, the de Sitter spacetime is characterized by the fact that its spatial Killing horizons have minimal modes for the Dirac operator. As a consequence, the de Sitter spacetime is the only vacuum of this type for which the induced metric tensor on some of its Killing horizons is at least equal to that of a round $(n-1)$-sphere. This extends uniqueness theorems shown in \cite{BGH,C} by Boucher-Gibbons-Horowitz and 
Chru{\'s}ciel to more general horizon metrics and to the non-single horizon case.     
\end{abstract}



\thanks{The second author was partially 
supported by a Spanish MEC-FEDER grant No. MTM2011-22547}

\date{\today}   

\maketitle \pagenumbering{arabic}
 
\section{Introduction}
An $(n+1)$-dimensional {\em vacuum} spacetime with  cosmological constant $\Lambda$ is a Lorentzian manifold $({\mathcal V},
g_{ab})$ satisfying the Einstein equation $R_{ab}=\Lambda\,{g}_{ab}$, where $R_{ab}$ is the Ricci tensor of
the metric $g_{ab}$. The vacuum is said to be
{\em static} when  
\begin{equation}\label{static}
{\mathcal V}=\R\times M,\qquad ds^2=-V^2\,dt^2+{^n\!}g_{ab}dx^adx^b,
\end{equation}
where $(M,{^n\!}g_{ab})$ is an $n$-dimensional connected Riemannian manifold, that we will take to be orientable, 
standing for the 
unchanging slices of constant time and $V\in C^\infty (M)$ is a  non-trivial smooth function
 on $M$. In the case  of positive cosmological constant $\Lambda>0$, it seems physically natural to
 require   spatially
 compact solutions, that is, the Cauchy hypersurface $M$ is usually 
 taken to be {\em compact}. The vacuum Einstein equation can be translated into the following two conditions on $(M,{^n\!}g_{ab})$ and $V$: 
\begin{equation}\label{vacuum}
\nabla_a\nabla_bV=V\left({^n\!}R_{ab}-\Lambda\, {^n\!}g_{ab}\right),\qquad \nabla^2 V=-\Lambda\, V,
\end{equation}
where ${^n\!}R_{ab}$ is the Ricci tensor, $\nabla_a$ is the covariant derivative  and \hfill\break $\nabla^2=
\nabla_a\nabla^a={^n\!}g^{ab}\nabla_a\nabla_b$
is the Laplace operator of the Riemannian manifold $(M,{^n\!}g_{ab})$.  By the second equation, we see that the 
spacetime cannot be globally
static, that is, the lapse function $V$ 
changes sign on $M$.
This means that the causal Killing field $\frac{\partial}{\partial t}$ is
lightlike on some locus (in fact, on a hypersurface, as we will point out below) of ${\mathcal V}$. Taking traces in
the first  equation in (\ref{vacuum}) and taking into account the second one, we conclude
immediately that 
\begin{equation}\label{scalar}
{^n\!}R=(n-1)\Lambda,
\end{equation} 
where ${^n\!}R$ is the scalar curvature of $(M,{^n\!}g_{ab})$.\\

Denote by $\Sigma\subset M$ the zero set of the lapse function $V$. If there exists a point 
$x\in\Sigma$ with $(\nabla_a V)(x)=0$, then the unicity of the solutions to the first (integrable linear) 
equation in (\ref{vacuum}) would imply that $V$ is identically zero on $M$. This means that 
$\Sigma$ is a closed smooth hypersurface embedded in $M$, although it may have several connected 
components. Then the null hypersurface ${\mathbb R}\times \Sigma$ of the spacetime ${\mathcal V}$ is called a {\em Killing horizon},  because
it is the locus of 
${\mathcal V}$ where the causal non vanishing Killing vector field $\frac{\partial}{\partial t}$  is lightlike. 
From a physical point of view, one usually says that ${\mathbb R}\times\Sigma$ is a {\em cosmological event horizon}, 
the horizon occurring because of the rapid expansion of the space due to the  $\Lambda$ term. {\em A common abuse of language will allow us to call {\rm horizon} to the projected hypersurface $\Sigma$ of $M$}. Indeed, each component $\Omega$ of $M$ where $V$ is positive stands for a spatial region
where communication is possible and the components $\Sigma_\alpha$, $\alpha=1,\dots,k$, of its boundary $\partial \Omega=\Sigma$ represent unattainable barriers for signals. 
From now on, we will restrict ourselves to work on open domains $\Omega\subset M$ of this type with
compact closure $\overline{\Omega}$ and non empty boundary $\Sigma=\partial \Omega$, not 
necessarily connected. Since 
the Riemannian manifold with non empty boundary $(\Omega,{^n\!}g_{ab})$, with ${^n\!}g_{ab}\in
C^\infty(\overline{\Omega})$ and $V\in
C^\infty(\overline{\Omega})$, positive on $\Omega$ and vanishing along $\Sigma$, completely determine the 
physically realistic regions of the vacuum spacetime $({\mathcal V},g_{ab})$. It is also usual to call the triple
$(\Omega,{^n\!}g_{ab},V)$ a {\em positive static triple}. \\ 

The paradigmatic example of a positive static triple with cosmological constant $\Lambda>0$ is given by choosing $(\Omega,{^n\!}g_{ab})=(\S_+^n
\left(\sqrt{\frac{n}{\Lambda}}\right),{^n\!}\delta_{ab})$ the open upper $n$-hemisphere in $\R^{n+1}$ of radius $\sqrt{\frac{n}{\Lambda}}$, where ${^n\!}\delta_{ab}$ is the Euclidean metric tensor. In particular, ${^n\!}R_{ab}=\frac{\Lambda}{n}(n-1){^n\!}\delta_{ab}$ and $\Sigma=\partial\Omega={\Bbb S}^{n-1}\left(\sqrt{\frac{n}{\Lambda}}\right)$ is the equator. It is
easy to see that, if ${\bf p}$ is the pole of ${\mathbb S}_+^n\left(\sqrt{\frac{n}{\Lambda}}\right)$,  the corresponding {\em height} function $h$ given by$$
h(x)=x^a {\bf p}_a,\qquad \forall x\in{\mathbb S}_+^n\left(\sqrt{\frac{n}{\Lambda}}\right)$$
is positive on ${\mathbb S}_+^n\left(\sqrt{\frac{n}{\Lambda}}\right)$, vanishes along $\Sigma= {\Bbb S}^{n-1}\left(\sqrt{\frac{n}{\Lambda}}\right)$ and satisfies the Obata type equation (see \cite{K,O})$$
\nabla_a\nabla_b h=-\frac{\Lambda}{n}h\;{^n\!}\delta_{ab}.$$
As a consequence, $\nabla^2 h=-\Lambda\,h$, and so equations (\ref{vacuum}) are satisfied. Indeed, the corresponding spacetime $$(dS_+,ds^2)= (\R\times\S_+^n,-h^2\,dt^2+
{^n\!}\delta_{ab}
dx^adx^b)$$ is nothing but an open domain in
de Sitter spacetime of radius $\sqrt{\frac{n}{\Lambda}}$ bounded by a certain lightlike hypersurface. \\

In fact, the so called {\em cosmic no-hair conjecture}, formulated
by Boucher-Gibbons-Horowitz (see \cite[p.\,2449]{BGH}) refers to a postulated and de\-si\-red uniqueness 
for the above example: 

\begin{quote}
{\em The only $(n+1)$-dimensional static vacuum spacetime with  $\Lambda>0$ and connected cosmic event horizon is the
de Sitter spacetime of radius $\sqrt{\frac{n}{\Lambda}}$. In other words, the only $n$-dimensional positive static triple $(\Omega,{^n\!}g_{ab},V)$ with single-horizon $\Sigma=\partial\Omega$ and cosmological constant $\Lambda>0$ is 
given by a round hemisphere $(\mathbb S_+^n\left(\sqrt{\frac{n}{\Lambda}}\right),{^n\!}\delta_{ab})$ of radius $\sqrt{\frac{n}{\Lambda}}$, where the lapse function $V$ is taken as the height function attaining its maximum at the pole.}
\end{quote}

Connectedness of the horizon $\Sigma$ is essential for this conjecture to be true. In fact, we can easily construct positive vacuum 
triples $(\Omega,{^n\!}g_{ab},V)$ with double-horizon $\Sigma$. It suffices to take $\Omega=
\mathbb S_+^1\left(\frac1{\sqrt{\Lambda}}\right)\times P$, where $\mathbb S_+^1(\frac1{\sqrt{\Lambda}})$ is the upper half-circle
of radius $\frac1{\sqrt{\Lambda}}$ and $P$ is any \hfill\break $(n-1)$-dimensional Einstein compact manifold with Ricci curvature equal 
to $\Lambda$ (for example, the choice of $P$ as the sphere ${\mathbb S}^{n-1}\left(\sqrt{\frac{n-2}{\Lambda}}\right)$ provides
the so called Nariai spacetime), 
${^n\!}g_{ab}$ the product metric and  $V(t,x)=\sin \sqrt{\Lambda}\,t$
for all $\frac1{\sqrt{\Lambda}}e^{i\sqrt{\Lambda}\,t}\in  \mathbb S_+^1\left(\frac1{\sqrt{\Lambda}}\right)$ and all $x\in P$.
One can check (see, for example, \cite[p.\,51, 2.]{La2}) that these choices ensure that equations (\ref{vacuum})
are satisfied. Moreover, it
is immediate that, in this case, $\Sigma$ consists just of two copies of $P$.\\

The cosmic no-hair conjecture is closely related to another conjecture by Fischer and Marsden concerning the Riemannian
metrics which are critical points for the scalar curvature map (see \cite{B,FM,QY}). Indeed, we may rephrase
this Fischer-Marsden conjecture (see \cite[{Conjecture 2.}]{Sh}) just as  the  Boucher-Gibbons-Horowitz conjecture above by 
removing the single-horizon assumption. So, the last  aforementioned examples should be viewed as counterexamples 
of this conjecture, at least as long as stated in such a generality. \\

Since the gradual spreading of the cosmic no-hair conjecture, some results about the subject have been interpreted as advances lending support to it. Firstly, it
is clear from (\ref{vacuum}) that, if $(\Omega,{^n\!}g_{ab})$ is Einstein (or even if just the Ricci tensor ${^n\!}R_{ab}$ has $\Lambda$ as 
a lower or upper bound), the Obata type theorem in \cite{Re} (cf.\,\cite{O}) solves the conjecture in the
affirmative, even in the case where the horizon $\Sigma$ is not assumed in principle to be connected. The same positive answer is obtained when  $(\Omega,{^n\!}g_{ab})$ is supposed to be conformally flat, a result 
proved independently by Kobayashi in \cite{Ko} and Lafontaine in \cite{La1}. Moreover, Chru{\'s}ciel, generalizing some early computations by Lindblom in \cite{Li}, showed in \cite{C} the following integral inequality 
\begin{equation}\label{scalar-integral}
\sum_{\alpha=1}^k\kappa_\alpha\int_{\Sigma_\alpha}\left({^{(n-1)}\!}{R_\alpha}- \frac{\Lambda}{n}(n-1)(n-2)\right)\,d\Sigma_\alpha\ge 0,
\end{equation} 
where ${^{(n-1)}\!}{R_\alpha}$ is the scalar curvature of the metric tensor ${^{(n-1)}\!}h_{ij}^\alpha$ induced from $\Omega$ on the connected component $\Sigma_\alpha$ of the event horizon
$\Sigma$ and each constant $\kappa_\alpha>0$ is the corresponding {\em surface gravity} 
on $\Sigma_\alpha$ (see (\ref{surface-gravity}) for a definition). The equality implies that   $\Omega$
is the  round hemisphere of radius $\sqrt{\frac{n}{\Lambda}}$ (and we have $k=1$ a fortiori, that is, $\Sigma$ must be connected).\\ 

Inequality (\ref{scalar-integral}) has two important consequences. Clearly, it implies that at least one of the integrals 
in the sum must be non negative. Since an $(n-1)$-dimensional round sphere of radius $\sqrt{\frac{n}{\Lambda}}$ has
constant scalar curvature $\frac{\Lambda}{n}(n-1)(n-2)$, at least one of the components $\Sigma_\alpha$ of the event horizon $\Sigma$
has total scalar curvature greater than or equal to that of such a sphere. In the single-horizon case, this means that: 
  
\begin{quote}
{\em The de Sitter vacuum minimizes the integral of the scalar curvature of the induced metric on the event horizon
among all the single-horizon positive static triples with the same cosmological constant.}
\end{quote}

 In the case $n=3$, inequality (\ref{scalar-integral}) along with the considerations above imply that the horizon has at least one
 genus zero component and, in the single-horizon case, one obtains the inequality
\begin{equation}\label{area}
12\,\pi \ge \Lambda\,\hbox{\rm area\,}(\Sigma),
\end{equation}
discovered by Boucher-Gibbons-Horowitz in \cite{BGH} (see also \cite{Sh} for an analogous computation in the case of
multiple horizons). As a consequence of these last inequalities, one has the following uniqueness result by Chru{\' s}ciel (see
\cite{C,Sh}):

\begin{quote} {\em 
Let $(\Omega,{^n\!}g_{ab},V)$ be an $n$-dimensional positive static triple  with cosmological constant $\Lambda>0$ and suppose that the 
horizon $(\Sigma,{^{(n-1)}\!}h_{ij})$   is isometric to a sphere $\mathbb S^{n-1}\left(r\right)$
of radius $r>0$. Then $r\le \sqrt{\frac{n}{\Lambda}}$, and if the equality holds, the triple $(\Omega,{^n\!}g_{ab},V)$ is   
given by a round hemisphere $(\mathbb S_+^n\left(\sqrt{\frac{n}{\Lambda}}\right),{^n\!}\delta_{ab})$ of radius $r=\sqrt{\frac{n}{\Lambda}}$, where the lapse function $V$ is  a height function attaining its maximum at the pole.
}\end{quote}

In fact, in order to have the  uniqueness result above, it would be enough that the scalar curvature ${^{(n-1)}\!}
{R}$ of $(\Sigma,{^{(n-1)}\!}h_{ij})$ is at most 
equal to that of the sphere  $\mathbb S^{n-1}\left(\sqrt{\frac{n}{\Lambda}}\right)$, that is, 
${^{(n-1)}\!}
{R}\le \frac{(n-1)(n-2)}{n} \Lambda$.\\

Eventhough all these results have been thought of as evidences for solving affirmatively the cosmic no-hair conjecture, 
 at least in the single-horizon case, they should be viewed as signs that the desired uniqueness of the de Sitter
 spacetime seems to require some boundary condition, that is, some topological, geometrical or
 analytical assumptions on the cosmic event horizon, as in the case of negative cosmological constant $\Lambda$
 (see \cite{BGH,An0,AD,CH,CS,GSW,HM1,Le,Wa1,Wa2}) or zero (see \cite{An,BM,C1,I}). This point of view has been confirmed by
 Gibbons-Hartnoll-Pope, who constructed in \cite{GHP} counterexamples to the 
cosmic no-hair conjecture, in the cases $4\le n\le 8$, by using the Riemannian inhomogeneous Einstein metrics
found by B{\"o}hm in \cite{Bo} (as well as by the uniqueness resuts obtained in \cite{BGH,C,LR,Ma}). In these counterexamples, one can find event horizons which are topologically
spherical but endowed with non round metrics and Riemannian products of spheres. Anyway, Gibbons-Hartnoll-Pope 
showed that these examples are dinamically unstable and should evolve into an asymptotically de Sitter spacetime.\\

Therefore, taking into account the remarks above and the
uniqueness results obtained by Chru{\'s}ciel in \cite{C} and by Boucher-Gibbons-Horowitz in \cite{BGH}, 
it becomes interesting to find necessary conditions, like inequalities (\ref{scalar-integral}) and (\ref{area}), on a given $(n-1)$-dimensional compact Riemannian manifold
$(\Sigma,{^{(n-1)}\!}h_{ij})$ in order to be a connected component of the event horizon of a positive static triple $(\Omega,{^n\!}g_{ab},V)$ with
cosmological constant $\Lambda>0$. In this way, we will be able to approach to the uniqueness of the Sitter spacetime imposing 
natural conditions on the event horizon and ruling out the single-horizon assumption.\\

Indeed, in this paper, in the first place we will  prove that the de Sitter spacetime minimizes the absolute value of the modes of the Dirac operator
on each component of the event horizon among all the positive static triples. Precisely, we will show the following result.

\begin{theorem}\label{main-th} 
Let $(\Omega,{^n\!}g_{ab},V)$ be
a positive static  triple with cosmological constant $\Lambda>0$. Suppose
that $\Omega$ is a spin manifold (this is always the case if $n=3$) and that  ${\partial\!\!\!/}_\alpha$ and $\kappa_\alpha$ are respectively the Dirac operator 
of the Riemannian spin structure induced on the connected component $\Sigma_\alpha$, $\alpha=1,\dots,k$, of the 
event horizon $\Sigma=\partial\Omega$ and the surface gravity of $\Sigma_\alpha$. Then  $$
\left|\lambda ({\partial\!\!\!/}_\alpha)\right|\ge \frac{n-1}{2}\sqrt{\frac{\Lambda}{n}}
\left(\frac{\kappa_\alpha}{\kappa_{\rm max}}\right),
$$
for all the modes $\lambda$ of ${\partial\!\!\!/}_\alpha$, where $\kappa_{\rm max}=
\max_{\alpha=1,\dots,k}\kappa_\alpha$. If equality holds for some $\alpha=1,\dots,k$, then $k=1$, that is, the event horizon $\Sigma$ is connected, $\kappa_\alpha=\kappa_{\rm max}$, the Riemannian manifold $(\Omega,{^n\!}g_{ab})$ is   
given by a round hemisphere $\left(\mathbb S_+^n\left(\sqrt{\frac{n}{\Lambda}}\right),{^n\!}\delta_{ab}\right)$ of radius $\sqrt{\frac{n}{\Lambda}}$ and the lapse function $V$ is a height function attaining its maximum at the pole. 
\end{theorem}

Note that the lower bound obtained in Theorem \ref{main-th} is independent of scale changes
 in the lapse function $V$. In fact, it is obvious from (\ref{vacuum}) that, if $V$ is a lapse function, any
multiple of $V$ is a lapse function as well. \\ 

As a consequence, and by using a nice upper estimate of the Dirac operator
on spheres due to Herzlich (see \cite{He}), who improved an original general estimate by Vafa-Witten (see \cite{VW}), we obtain another result which generalises the 
Chru{\'s}ciel and Boucher-Gibbons-Horowitz theorem in two directions: it applies to the general case of non single-horizon $\Sigma$ and moreover there is no need to impose that the metric tensor on the horizon is the round spherical metric, but only that it dominates this round metric. 

\begin{theorem}\label{uniqueness}
Let $(\Omega,{^n\!}g_{ab},V)$ be an $n$-dimensional positive static triple  with cosmological constant $\Lambda>0$ and suppose that 
$\Omega$ is spin and that there is a connected component $\Sigma_0$ of the event horizon $\Sigma=\partial\Omega$ which  is diffeomorphic to an $(n-1)$-dimensional sphere. If the 
corresponding induced tensor metric ${^{(n-1)}\!}h_{ij}^0$ is pointwise at least equal to that of the round standard sphere of radius $r>0$, then $r\le \sqrt{\frac{n}{\Lambda}}\left(\frac{\kappa_{\rm max}}{\kappa_0}\right)$. If equality is attained, then $\Sigma$ is connected, $\kappa_0=\kappa_{\rm max}$, the Riemannian manifold $(\Omega,{^n\!}g_{ab})$ is   
given by a round hemisphere $\left(\mathbb S_+^n\left(\sqrt{\frac{n}{\Lambda}}\right),{^n\!}\delta_{ab}\right)$ of radius $r=\sqrt{\frac{n}{\Lambda}}$ and the lapse function $V$ is a height function attaining its maximum at the pole.
\end{theorem}

\section{Geometry of a modification of the Fermat conformal metric}
Given a positive static triple $(\Omega,{^n\!}g_{ab},V)$, the conformal metric \hfill\break ${\widetilde {^n\!g}_{ab}}=\frac1{V^2}{^n\!} g_{ab}$
on $\Omega$ is known as the {\em Fermat (or optical) metric}. This is why the geodesics of this new metric are the spatial projections 
of the light rays in the corresponding spacetime $(\mathbb R\times \Omega, -V^2\,dt^2+{^n\!}g_{ab}dx^adx^b)$. The geometrical features 
of this conformal metric have been explicitly or implicitly studied in order to analyse the behaviour of static spacetimes, mainly with non-null
cosmological constant (see, for example, \cite{BGH,CS,GSW,Li,Q}). In \cite{C,HM1}, a suitable modification of this metric has been used. We will
consider again this modified Fermat metric in order to prove Theorems \ref{main-th} and \ref{uniqueness}.

\begin{lemma}\label{R^*}
Let $(\Omega,{^n\!}g_{ab},V)$ be a positive static vacuum and $\varepsilon> 0$ an arbitrary positive real number. Then the conformal Riemannian metric ${^n\!}g_{ab}^*=\left(\frac{\varepsilon}{V+\varepsilon}\right)^{2}{^n\!}g_{ab}$
has scalar curvature
\begin{equation}\label{estimate*}
{^n\!}R^*=\frac{(n-1)\Lambda}{\varepsilon^2}\left(\varepsilon^2-V^2-\frac{n}{\Lambda}W\right),
\end{equation}
where $W$ is the squared norm of the gradient of the lapse function $V$, that is, $W={^n\!}g^{ab}\nabla_a V\nabla_b V$.
Moreover, ${^n\!}g_{ab}^*$ coincides with ${^n\!}g_{ab}$ on the event horizon $\Sigma=\partial\Omega$.
\end{lemma}

\pf
The last assertion is clear because $V$ vanishes along $\Sigma$. If we rewrite the conformal change between the metrics ${^n\!}g_{ab}$ and ${^n\!}g_{ab}^*$ as $$
{^n\!}g_{ab}=\left(\frac{V+\varepsilon}{\varepsilon}\right)^2{^n\!}g_{ab}^*,$$
then the relation between the corresponding Ricci tensors ${^n\!}R_{ab}^*$ and ${^n\!}R_{ab}$
 on the compact manifold ${\overline\Omega}$, is given by (see \cite[p.\! 59]{Be}):
$$ 
 {^n\!}R_{ab}^*={^n\!}R_{ab}+(n-2)\frac{\nabla_a\nabla_b V}{V+\varepsilon} + 
 \frac{\nabla^2 V}{V+\varepsilon}{^n\!}g_{ab}-(n-1)\frac{W}{(V+\varepsilon)^2}{^n\!}g_{ab}.
$$
Taking traces with respect to ${^n\!}g_{ab}$ and multiplying by $(V+\varepsilon)^2$, we obtain
the corresponding relation for their scalar curvatures$$
\varepsilon^2({^n\!}R^*)={^n\!}R(V+\varepsilon)^2+2(n-1)(V+\varepsilon)\nabla^2 V-n(n-1)W.
$$
As we pointed out in (\ref{vacuum}) and (\ref{scalar}), since ${^n\!}R=(n-1)\Lambda$ and
the function $V$ is an eigenfunction of $\nabla^2$ associated with the
eigenvalue $-\Lambda$, we finally
get the required expression for the scalar curvature ${^n\!}R^*$ of the 
compact Riemannian manifold with boundary $({\overline\Omega},{^n\!}g_{ab}^*)$. 
\qed\\

Our aim now is to prove that the conformal metric ${^n\!}g_{ab}^*$ has non-negative scalar curvature 
${^n\!}R^*$. The proof of this assertion will be a reformulation of an identity due to Lindblom (\cite{Li})
showing that ${^n\!}R^*$ satisfies a certain elliptic second order equation on ${\overline\Omega}$ 
and a short analysis displaying the geometry of the event horizon $\Sigma$ as a 
hypersurface of ${\overline\Omega}$. Let us start by this last point. In fact, in the discussion after equality 
(\ref{scalar}), we noticed that along the event horizon $\Sigma$, which is the zero set of the lapse function $V$,
the gradient vector field $\nabla^a V$ has no zeroes and so, after normalization, it provides an inner normal vector field
$N^a$ for the hypersurface $\Sigma$ in ${\overline\Omega}$. Then, it is well-known that the extrinsic curvature
(or second fundamental form) $p_{ij}=-{^{n}\!}g_{ic}\nabla_jN^c$ of  any regular level set of $V$ is given by the equation$$
\nabla_i\nabla_j V+(N^c\nabla_c V)p_{ij}=0.$$
Since $N^c\nabla_c V$ has no zeroes 
on $\Sigma$ (which is the regular level set $V^{-1}(\{0\})$),
by the first equality in (\ref{vacuum}), we have 
\begin{equation}\label{totally-geodesic}
\hbox{$p_{ij}=0$ on $\Sigma$, that is, {\em
$\Sigma$ is totally geodesic in $\overline\Omega$.}}
\end{equation} 
Then, the fact that $p_{ij}=0$ and (\ref{vacuum}) give$$
\nabla_i(N^b\nabla_b V)_{|\Sigma}= -{^{n}\!}g^{bj}p_{ij}\nabla_b V+N^b(\nabla_i\nabla_b V)_{|\Sigma}=0,$$
and so the normal derivative $N^a\nabla_a V$ is constant on each connected component of the
event horizon $\Sigma$. Recall that  $\Sigma_\alpha$, $\alpha=1,\dots,k$, denote its  connected components. 
So, we dzfinz
\begin{equation}\label{surface-gravity}
\kappa_\alpha\equiv (N^a\nabla_a V)_{|{\Sigma_\alpha}}>0,\qquad\alpha=1,\dots,k.
\end{equation}
Note that each $\kappa_\alpha$ is positive because, as we pointed out in the discussion following 
equality (\ref{scalar}), the gradient of $V$ never
vanishes on its zero set. Physically,
each $\kappa_\alpha$ is the {\em surface gravity} on the component $\Sigma_\alpha$ of the event horizon
of the $n$-dimensional spatial slice $\Omega$.

\begin{proposition}\label{compactification}
Let $(\Omega,{^n\!}g_{ab},V)$ be an $n$-dimensional positive static triple with cosmological
constant $\Lambda>0$ and let $\Sigma_\alpha$, $\alpha=1,\dots,k$, be the connected components of
its event horizon.  Then, taking \hfill\break$\varepsilon =
\sqrt{\frac{n}{\Lambda}}\kappa_{\rm max}>0$, where $\kappa_{\rm max}$ is
the greatest surface gravity among all the components $\Sigma_\alpha$,  yields $$
{^{(n-1)}\!}h_{ij}^*={^{(n-1)}\!}h_{ij},\quad {^n\!}R^*\ge 0,\quad (p^*)_{|\Sigma_\alpha}
=(n-1)\sqrt{\frac{\Lambda}{n}}\left(\frac{\kappa_\alpha}{\kappa_{\rm max}}\right),$$
where  ${^n\!}R^*$ is the scalar curvature of the conformal metric ${^n\!}g_{ab}^*$ defined in Lemma \ref{R^*}, $p^*$ is the trace of the (inner) extrinsic curvature  $p_{ij}^*$ of the event horizon $\Sigma=\partial \Omega$ as a hypersurface
of $({\overline\Omega},{^n\!}g_{ab}^*)$. 
\end{proposition}

\pf 
 The first equality is simply the last assertion in the statement of Lemma \ref{R^*}. As for the second one, we define
a function $\Phi$  on
${\overline\Omega}$ by$$
\Phi=\varepsilon^2-V^2-\frac{n}{\Lambda}W=\varepsilon^2-V^2-\frac{n}{\Lambda}\;{^n\!}g^{ab}\;\nabla_a V\nabla_b V,
$$
 that is, up to a positive constant, the
righthand side of (\ref{estimate*}). Taking into account the more or less explicit computations 
in \cite{Li,C,Q} or the classical Bochner formula for the Laplace operator of
the squared length of a gradient, we have
\begin{eqnarray*}
\frac1{2}\nabla^2 \Phi & = & -W-V\nabla^2 V-\frac{n}{\Lambda}(\nabla_a\nabla_b V)(\nabla^a\nabla^b V) \\ &-&
\frac{n}{\Lambda}\,\nabla_a V\nabla^a(\nabla^2 V)-\frac{n}{\Lambda}\,{^n\!}R_{ab}\,\nabla^a V\nabla^b V.
\end{eqnarray*} 
Using (\ref{vacuum}), we obtain
\begin{eqnarray}\label{Bochner}
\frac1{2}\nabla^2 \Phi&=& -\frac{n}{\Lambda}\left(\nabla_a\nabla_b V+\frac{\Lambda}{n}V\,{^n\!}g_{ab}\right)\left(
\nabla^a\nabla^b V+\frac{\Lambda}{n}V\,{^n\!}g^{ab}\right)\\&-&
\frac{n}{\Lambda}\frac{(\nabla_a\nabla_bV ) (\nabla^a V\nabla^b V)}{V}-W.\nonumber
\end{eqnarray}
On the other hand, from the  definition of $\Phi$, we have
\begin{equation}\label{gradient-Fi}
\nabla_a\Phi\nabla^a V=-2VW-2\frac{n}{\Lambda}\;\nabla_a \nabla_b V\,\nabla^a V\nabla^b V.
\end{equation}
Hence, by combining (\ref{Bochner}) and (\ref{gradient-Fi}), we finally obtain
\begin{eqnarray}\label{Bochner2}
\frac1{2}\nabla^2 \Phi&=& -\frac{n}{\Lambda}\left(\nabla_a\nabla_b V+\frac{\Lambda}{n}V\,{^n\!}g_{ab}\right)\left(
\nabla^a\nabla^b V+\frac{\Lambda}{n}V\,{^n\!}g^{ab}\right)\\&+&
\frac1{2V}\nabla_a\Phi\nabla^aV.\nonumber
\end{eqnarray}
A standard application of min-max principle to this elliptic equation (\ref{Bochner2}) implies that 
$\Phi$ reaches its minimum
value at $\Sigma$.    
Then, to control the  behaviour of the function $\Phi$  along the event horizon, it suffices to recall
that $\Sigma=V^{-1}(\{0\})$ and that$$
W_{|\Sigma}=(\nabla_aV\nabla^aV)_{|\Sigma}= (N^a\nabla_aV)^2\le 
\frac{\Lambda}{n}\varepsilon^2,$$
according to (\ref{surface-gravity}) and the choice of $\varepsilon$. Therefore, we conclude that$$
\Phi_{|\Sigma}=\left(\varepsilon^2-V^2-\frac{n}{\Lambda}W\right)_{|\Sigma}\ge 0,$$
and that, if $\Phi$ is not constant,  the equality holds only on the connected components of $\Sigma$ with maximum
surface gravity. This proves that $\Phi\ge 0$ on the whole of $\Omega$ and so, using (\ref{estimate*}) in Lemma
\ref{R^*}, we obtain the second assertion ${^n\!}R^*\ge 0$.\\

To finish the proof, it remains to compute the mean curvature $p^*$ of the event
horizon $\Sigma$ as a hypersurface of the compact Riemannian manifold
$({\overline\Omega},{^n\!}g_{ab}^*)$. Observe that, by definition in Lemma \ref{R^*}: $$
{^n\!}g_{ab}^*=\left(\frac{\varepsilon}{V+\varepsilon}\right)^2{^n\!}g_{ab}.$$
But,  we know from (\ref{totally-geodesic}) that  the extrinsic curvature $p_{ab}$ of the event horizon 
$\Sigma$ as a hypersurface of $({\overline \Omega},{^n\!g}_{ab})$ vanishes (in fact, $\Sigma$ is
 a totally geodesic hypersurface). 
So, in order to compute $p^*$, it suffices to use the well-known relation between the two mean curvatures 
of a hypersurface corresponding to two metrics on the ambient space which are conformal (see, for instance, \cite{E}):$$
p^*=\frac1{\rho}\big(p- (n-1)N^a\nabla_a\log \rho\big)=-\frac{n-1}{\rho^2}
N^a\nabla_a\rho,$$
where $\rho=\frac{\varepsilon}{V+\varepsilon}$ and $N^a$ is the inner unit normal  along $\Sigma$ with respect to the
metric ${^n\!}g_{ab}$. Since we have$$
(\nabla_a \rho)_{|\Sigma}=-\left.\frac{\varepsilon \nabla_a V}
{\big(V+\varepsilon \big)^2}\right|_{\Sigma}=-\frac{\rho^2}{\varepsilon}(\nabla_a V)_{|\Sigma}.$$  
Thus, we obtain$$
p^*=\frac{n-1}
{\varepsilon}(N^a\nabla_a V)_{|\Sigma}.$$   
Taking into account (\ref{surface-gravity}) and the choice of $\varepsilon$, for any  $\alpha=1,\dots,k$, we finally have$$
(p^*)_{|\Sigma_\alpha}=\frac{n-1}
{\varepsilon}\kappa_\alpha=(n-1)\sqrt{\frac{\Lambda}{n}}\left(\frac{\kappa_\alpha}{\kappa_{\rm max}}
\right).$$
\qed

\begin{remark}{\rm
Note that equation (\ref{Bochner2}) is equivalent to$$
\frac1{2}\hbox{\,\rm div\,} \frac1{V}\nabla \Phi=-\frac{n}{\Lambda}\left(\nabla_a\nabla_b V+\frac{\Lambda}{n}V\,{^n\!}g_{ab}\right)\left(
\nabla^a\nabla^b V+\frac{\Lambda}{n}V\,{^n\!}g^{ab}\right),$$
where it can easily checked that the vector field $\frac1{V}\nabla\Phi$ extends smoothly to the
boundary $\Sigma$.
So, we have the inequality$$
\hbox{\,\rm div\,} \frac1{V}\nabla \Phi\le 0.$$
By integrating and using the divergence theorem, we obtain$$
\int_\Sigma \left(\frac1{V}N^a\nabla_a\Phi\right)\, d\Sigma\ge 0,$$
where $N=\frac{\nabla V}{|\nabla V|}$ is the inner unit normal along $\Sigma$. Now, by using the 
definitions of $\Phi$ and of the $\kappa_\alpha$, $\alpha=1,\dots,k$, we have that$$
\left(\frac1{V}N^a\nabla_a\Phi\right)_{|\Sigma_\alpha}=\kappa_\alpha\left(
\frac{n}{\Lambda}{^n\!}R_{ab}N^aN^b-(n-1)\right).$$
Then, equations (\ref{scalar-integral}) and (\ref{area}) are direct consequences from the discussion
above and the Gau{\ss} equation$$
{^{(n-1)}\!}R=(n-1)\Lambda-2{^n\!}R_{ab}N^aN^b$$
on the totally geodesic hypersurface $\Sigma$.  
}
\end{remark}

\section{Proofs of Theorems \ref{main-th} and \ref{uniqueness}}
Suppose now that the positive static triple $(\Omega, {^n\!}g_{ab},V)$ is such that the compact orientable 
$n$-dimensional manifold $\overline{\Omega}$ with non empty boundary is a spin manifold (this is always the case when the spatial slice $M$ of the spacetime is spin) and  that
we have fixed a spin structure on it.  Since the horizon 
$\Sigma=\partial\Omega $ is always an orientable hypersurface, it is also a 
spin manifold and   that an induced spin structure on the horizon is inherited from the structure fixed on
${\overline \Omega}$. Moreover, for the Riemanian metric ${^n\!}g_{ab}$ on ${\overline \Omega}$ we have 
an associated spinor bundle $({\mathbb S}{\overline \Omega},\gamma^a,{\overline\nabla}_a,{\overline\partial})$,
where $\gamma^a$ are the Pauli matrices, ${\overline\nabla}_a$ the covariant derivative and ${\overline\partial}$ is the corresponding Dirac operator (for generalities on spin structures see any 
of \cite{BFGK,Fr,LM,BHMMM}). It is a well-known
fact that the restriction of the spinor bundle ${\mathbb S}{\overline \Omega}$  to the hypersurface $\Sigma$ 
can be identified with one or two copies of the spinor bundle corresponding to the induced spin structure and 
the induced Riemannian metric ${^{(n-1)}\!}h_{ij}$ according to the parity of the dimension $n$
of ${\overline \Omega}$. More precisely,  we have an isomorphism$$
\left({\mathbb S}{\overline\Omega}_{|\Sigma},\gamma^n{\overline\partial}+N^a
{\overline\nabla}_a-\frac1{2}p\right)\cong
\left\{\begin{array}{ll}
({\mathbb S}\Sigma,\partial\!\!\!/),\hbox{ if $n$ is odd} \\ \\
({\mathbb S}\Sigma,\partial\!\!\!/)\oplus({\mathbb S}\Sigma,
-\partial\!\!\!/),
\hbox{ if $n$ is even},
\end{array}\right.$$  
 where $N$ is the inner unit normal field along the horizon $\Sigma$, $p$ is the trace of its extrinsic curvature
 and $({\mathbb S}\Sigma,\partial\!\!\!/)$ are respectively the spinor bundle and the 
 Dirac operator corresponding to the 
 spin structure and to the Riemannian metric induced on $\Sigma$ (for this relationship between the spinor bundles on
 a hypersurface and on its ambient space, see, for instance, \cite{Ba3,Fr,HMZ1,HMZ2}). 
 Due to this identification we can say that each spinor field on 
 ${\overline \Omega}$ determines by restriction a spinor field on the event horizon $\Sigma$
 and we can talk about possible extensions to ${\overline \Omega}$ of the spinor fields defined
 on $\Sigma$. Furthermore,  from the identification between the operators $\partial
 \!\!\!/$ and $\gamma^n{\overline\partial}+N^a
{\overline\nabla}_a-\frac1{2}p$, it is immediate that, if $p$ is constant, the restriction to $\Sigma$ of a parallel spinor field on ${\overline\Omega}$
gives an eigenspinor on $\Sigma$ associated to the eigenvalue $-\frac1{2}p$ of the operator $\partial
\!\!\!/$. Taking into account the identifications above between bundles and operators, X. Zhang and the first two authors 
showed in \cite{HMZ1} (see also \cite[Theorem 3.7.1]{Gi}) that, if the scalar curvature ${^n\!}R$ of the metric 
${^n\!}g_{ab}$ is non negative on ${\overline\Omega}$ and the trace $p$ of the extrinsic curvature of $\Sigma$
in ${\overline\Omega}$ is also non negative, we have 
\begin{equation}\label{our-estimate}
|\lambda_1(\partial\!\!\!/)|\ge \frac{1}{2}\inf_\Sigma p,
\end{equation} 
where $\lambda_1(\partial\!\!\!/)$ stands for the eigenvalue of $\partial\!\!\!/$ with the lowest absolute value, 
and, if the equality holds, then the eigenspace
associated to $\lambda_1(\partial\!\!\!/)$ is built from  parallel spinor fields 
on ${\overline \Omega}$ (note that  in \cite{HMZ1,Gi} the inequality is given in terms of the normalized mean curvature of $\Sigma$).  It is straightforward to check that the approach  in \cite{HMZ1}
can be applied to each connected component $\Sigma_\alpha$, $\alpha=1,\dots,k$, of 
$\Sigma$ in the case where $\Sigma$ is not
connected. From 
equations (\ref{scalar}) and (\ref{totally-geodesic}), one can see that the estimate (\ref{our-estimate}) can be applied
to our situation. But, unfortunately, we get on each component, the obvious inequality $|\lambda_1(\partial\!\!\!/_\alpha)|\ge 0$. In fact, it is
clear that  (\ref{our-estimate}) is of interest 
only when $\inf_\Sigma p>0$. Nevertheless, we can obtain
some significant information by considering the conformal metric ${^n\!}g_{ab}^*$
on ${\overline\Omega}$ (see Lemma \ref{R^*} and Proposition \ref{compactification}). By combining the information provided
by Proposition \ref{compactification} and estimate (\ref{our-estimate}), we get$$
\left|\lambda_1 ({\partial\!\!\!/}_\alpha)\right|\ge \frac{n-1}{2}\sqrt{\frac{\Lambda}{n}}
\left(\frac{\kappa_\alpha}{\kappa_{\rm max}}\right),
$$
for each $\alpha=1,\dots,k$, as required in Theorem \ref{main-th}. As we mentioned before, according to 
\cite{HMZ1}, if the equality holds, then the eigenspace
associated with $\lambda_1(\partial\!\!\!/)$ is built of  parallel spinor fields 
on $({\overline \Omega}, {^n\!}g_{ab}^*)$. This implies the existence of  a non trivial  parallel spinor field
$\Psi\in \Gamma {\mathbb S}{\overline \Omega}$ with respect to the metric ${^n\!}g_{ab}^*$. 
It was shown by Hitchin in \cite{Hit} (see also \cite[Chapter 6]{BFGK}) that the existence
of a non-trivial parallel spinor forces the Ricci tensor to vanish everywhere. Then ${^n\!}R_{ab}^*
=0$ on ${\overline\Omega}$ and so ${^n\!}R^*=0$ as well. From (\ref{estimate*}) and 
(\ref{Bochner2}), we conclude 
\begin{equation}\label{Obata}
\nabla_a\nabla_b V+\frac{\Lambda}{n}\,V\;{^n\!}g_{ab}=0.
\end{equation}
Hence the compact Riemannian manifold with non empty boundary $({\overline\Omega},{^n\!}g_{ab})$ admits a
non trivial solution $V$ to the Obata  equation which is positive on $\Omega$ and vanishes on its boundary $\Sigma$.
Now, we may apply the  non empty boundary version of the Obata theorem due to Reilly (see \cite{Re}) and
conclude that  $({\overline\Omega},{^n\!}g_{ab})$ is a round hemisphere with radius $\sqrt{\frac{n}{\Lambda}}$
and $V$ is the height function with maximum at its center. This finishes the proof of Theorem \ref{main-th}.\\

As for Theorem \ref{uniqueness}, we suppose that there is a component $\Sigma_0$ of the event horizon $\Sigma$ which is diffeomorphic to
a sphere ${\mathbb S}^{n-1}$. Then the simple connectedness of this component implies that it supports only one spin structure and so the spin structure induced from ${\overline\Omega}$ 
on $\Sigma_0$ is the standard one on the $(n-1)$-dimensional sphere. On the other hand, we assume that
the metric ${^{(n-1)}\!}h_{ij}^0$ induced on $\Sigma$, after applying a diffeomorphism if necessary,  satisfies$$    
{^{(n-1)}\!}h_{ij}^0\ge {^{(n-1)}\!}\delta_{ij}, \qquad\hbox{pointwise},
$$
where ${^{(n-1)}\!}\delta_{ij}$ is the round metric of radius, say, $r>0$. This hypothesis
allows us to apply \cite[Theorem 1]{He} and conclude that$$  
\left|\lambda_1 ({\partial\!\!\!/}_0)\right|\le \frac{n-1}{2r},$$
where the equality implies that ${^{(n-1)}\!}h_{ij}^0$ is just the round metric. Combining this upper
bound for the eigenvalue of the Dirac operator ${\partial\!\!\!/}_0$ of the horizon with the least absolute value with the lower bound provided by Theorem \ref{main-th}, we have$$
\sqrt{\frac{\Lambda}{n}}
\left(\frac{\kappa_0}{\kappa_{\rm max}}\right)\le \frac1{r}.$$ 
If the equality is attained, we know that ${^{(n-1)}\!}h_{ij}^0= {^{(n-1)}\!}\delta_{ij}$ and, moreover, we may also apply the equality case 
in Theorem \ref{main-th}. So Theorem \ref{uniqueness} is proved.

\end{document}